\input amstex
\magnification=1200
\documentstyle{amsppt}
\NoRunningHeads
\NoBlackBoxes
\define\Mat{\operatorname{Mat}}
\define\so{\operatorname{\frak s\frak o}}
\define\gl{\operatorname{\frak g\frak l}}
\define\Hom{\operatorname{Hom}}
\define\xvr{\operatorname{\xi\varrho}}
\define\Der{\operatorname{Der}}
\define\ad{\operatorname{ad}}
\define\mYB{\operatorname{\widetilde mYB}}
\topmatter
\title On the nonHamiltonian interaction of two rotators. II.
Unravelling the algebraic structure
\endtitle
\author Denis V. Juriev
\endauthor
\affil\eightpoint\it Research Center for Mathematical Physics and
Informatics "Thalassa Aitheria"\linebreak
ul.Miklukho-Maklaya 20-180, Moscow 117437 Russia\linebreak
(e-mail: denis\@juriev.msk.ru)\linebreak
\endaffil
\endtopmatter
\document
This short note being a continuation of the first part [1] is devoted to the
unravelling of the algebraic structure, which governs the quadratic
nonHamiltonian interaction of two rotators described in [1]. It should
be considered in the context of a general ideology of the inverse problem
of representation theory [2]. The resulted objects are formalized as Lie
algebras with operators (Lie algebras with additional unary operations --
see [3]).

\subhead 1. Lie $\xvr$-algebras\endsubhead
First, let introduce and investigate the algebraic object, which has
the general algebraic meaning besides its relation to the discussed topic.

\definition{Definition 1} The {\it Lie $\xvr$-algebra\/} is a triple
$(\frak g,\xi,\varrho)$, where $\frak g$ is a Lie algebra, $\xi\in\Der(\frak
g)$, $\varrho$ is an operator in $\frak g$ commuting with $\xi$ such that
the following identities hold
$$\aligned
\varrho[x,y]_\varrho=[\varrho x,\varrho y],\quad
&\xi[x,y]_\varrho=[\varrho x,\xi y]+[\xi x,\varrho y],\\
\xi^2[x,y]_\varrho=[\varrho^2 x,y]&+[x,\varrho^2 y]-2[\varrho x,\varrho y]
\endaligned$$
where
$$[x,y]_\varrho=[\varrho x,y]+[x,\varrho y]-\varrho[x,y]-[\xi x,\xi y].$$
\enddefinition

\remark{Remark 1} The second identity of definition 1 may be rewritten as
$$\xi[\xi x,\xi y]=[Sx,y]+[x,Sy]-S[x,y],\qquad S=\xi\varrho.$$
\endremark

\remark{Remark 2} The bracket $[\cdot,\cdot]_\varrho$ coincides with the
the bracket $[\cdot,\cdot]_\Theta$ defined as
$$[x,y]_\Theta=\tfrac12\bigl([\Theta x,y]+[x,\Theta y]-\Theta[x,y]\bigr),
\qquad\Theta=2\varrho+\xi^2.$$
Moreover,
$$\Theta[x,y]_\Theta=[\varrho^2 x,y]+[x,\varrho^2y].$$
\endremark

\remark{Example 1} Let $\frak g=\gl(n)$ be the Lie algebra of all $n\times n$
matrices. Put $\xi_qx=qx-xq$, $\varrho_q x=qxq$ ($q\in\Mat_n$). The triple
$(\gl(n),\xi,\varrho)$ is a Lie $\xvr$-algebra.
\endremark

Note that the mapping $\varrho_q$ is quite familiar in the theory of
Jordan algebras and symmetric spaces [4,5].

\remark{Example 2} Let $*$ be an involution in $Mat_n$ and $\frak g=
\{x\in\Mat_n: x^*=-x\}$. Then the restriction of $\xi_q$ and $\varrho_q$
from example 1 onto $\frak g$ supplies it by the structure of a Lie
$\xvr$-algebra. In particular, $\so(n)$ is a Lie $\xvr$-algebra.
\endremark

\remark{Remark 3} $\xi$ is a derivation of the bracket
$[\cdot,\cdot]_\varrho$, i.e.
$$\xi[x,y]_\varrho=[\xi x,y]_\varrho+[x,\xi y]_\varrho.$$
\endremark

\remark{Exercise} To describe all structures of Lie $\xvr$-algebras on
$\so(3)$ [answer: each such structure coincides with one of example 2;
hint: use remark 3].
\endremark

This property is equivalent to the second identity of definition 1.

\proclaim{Theorem 1} In any $\xvr$-algebra the bracket
$[\cdot,\cdot]_\varrho$ is a Lie bracket compatible with $[\cdot,\cdot]$.
\endproclaim

Let us now relate the Lie $\xvr$-algebras to Lie bi-$\mYB$-algebras of
the article [6].

\definition{Definition 2 {\rm [6]}}

{\bf A.} The {\it Lie $\mYB$-algebra\/} is a Lie algebra $\frak g$ with the
bracket $[\cdot,\cdot]$, supplied by an operator $R:\frak g\mapsto\frak g$
such that
$$R[Rx,y]+R[x,Ry]=[Rx,Ry]+R^2[x,y].$$

{\bf B.} The {\it Lie bi-$\mYB$-algebra\/} is a Lie algebra $\frak g$ with
bracket $[\cdot,\cdot]$ supplied by two commuting operators $R_1$ and $R_2$
such that $(\frak g,R_1)$ and $(\frak g,R_2)$ are the Lie $\mYB$-algebras
with identical brackets $[\cdot,\cdot]_{R_1}$ and $[\cdot,\cdot]_{R_2}$.

{\bf C.} A Lie bi-$\mYB$-algebra $(\frak g,R_1,R_2)$ is called {\it
even-tempered} if the identities
$$\aligned
[R_1x,R_2y]+[R_2x,R_1y]-R_1R_2[x,y]&=[R_1^2x,y]+[x,R_1^2y]-R_1^2[x,y],\\
[R_1x,R_2y]+[R_2x,R_1y]-R_1R_2[x,y]&=[R_2^2x,y]+[x,R_2^2y]-R_2^2[x,y]
\endaligned
$$
hold.
\enddefinition

Examples of Lie $\mYB$-algebras and bi-$\mYB$-algebras were considered in
[6].

\proclaim{Proposition 1} Any even-tempered Lie bi-$\mYB$-algebra
$(\frak g,R_1,R_2)$ is a Lie $\xvr$-algebra with $\xi\!=\!R_1\!-\!R_2$,
$\varrho\!=\!R_1R_2$.
\endproclaim

The Lie $\xvr$-algebra of example 1 can be constructed from a Lie
bi-$\mYB$-algebras whereas one of example 2 can not.

\subhead 2. $\xvr$-structures on Lie algebras and quadratic
$I$-pairs\endsubhead In the definition of Lie $\xvr$-algebra it is natural
to consider $\xi$ as an inner derivation of $\frak g$. In this case the
defining identities of the Lie $\xvr$-algebra have the form
$$\aligned
\varrho[x,y]_\varrho=[\varrho x,\varrho y],\quad
&[q\,[x,y]_\varrho]=[\varrho x\,[q,y]]+[[q,x]\,\varrho y],\\
[q\,[q\,[x,y]_\varrho]]&=[\varrho^2 x,y]+[x,\varrho^2 y]-2[\varrho x,\varrho y]
\endaligned$$
where
$$[x,y]_\varrho=[\varrho x,y]+[x,\varrho y]-\varrho[x,y]-[[q,x]\,[q,y]]$$
and $\xi=\ad q$, $q\in\frak g$.

The following definition is self-explained.

\definition{Definition 3} The {\it $\xvr$-structure on the Lie algebra\/}
$\frak g$ is the set of operators $\varrho_q$ in $\frak g$ ($q\in\frak g$)
such that (1) $\varrho_q$ depends quadratically on $q$ (i.e.
$\varrho_{\lambda q}\!=\!\lambda^2\varrho_q$ and $\varrho_{q_1+q_2}\!+\!
\varrho_{q_1-q_2}\!=\!2\varrho_{q_1}\!+\!2\varrho_{q_2}$), (2) the triple
$(\frak g,\ad q,\varrho_q)$ is a Lie $\xvr$-algebra for all $q\in\frak g$.
\enddefinition

Below we shall consider $\frak g$-equivarint $\xvr$-structures on Lie
algebras $\frak g$.

\remark{Remark 4} The structures of Lie $\xvr$-algebras on $\gl(n)$ and
$\so(n)$ from examples 1,2 form $\xvr$-structures on them.
\endremark

\definition{Definition 4} A pair $(V_1,V_2)$ of linear spaces will be called
{\it I--pair\/} iff there are defined (nonlinear) mappings
$$\aligned &h_1:V_2\mapsto\Hom(\Lambda^2(V_1),V_1)\\
&h_2:V_1\mapsto\Hom(\Lambda^2(V_2),V_2)
\endaligned$$
such that $\forall A\in V_2$ its image $h_1(A)\in\Hom(\Lambda^2(V_1),V_1)$ is
a Lie bracket in $V_1$ and $\forall X\in V_1$ its image
$h_2(X)\in\Hom(\Lambda^2(V_2),V_2)$ is a Lie bracket in $V_2$. The Lie
bracket in $V_2$ corresponded to $A$ will be denoted by $[\cdot,\cdot]_A$,
whereas the Lie bracket in $V_1$ corresponded to $X$ will be denoted by
$[\cdot,\cdot]_X$.
\enddefinition

Here $\Hom(H_1,H_2)$ denotes the space of all linear operators from $H_1$ to
$H_2$, $\Lambda^2(H)$ is a skew square of the linear space $H$, so
$\Hom(\Lambda^2(H),H)$ is the space of all skew--symmetric bilinear binary
operations in $H$. The Lie brackets in $H$ form a submanifold
$\operatorname{Lie}(H)$ of the space $\Hom(\Lambda^2(H),H)$.

\proclaim{Theorem 2} Any $\xvr$-structure on the Lie algebra $\frak g$
induces a structure of the quadratic $I$-pair on $(\frak g,\frak g)$ with
$$[x,y]_q=[x,y]+[x,y]_{\varrho_q}.$$
\endproclaim

\remark{Remark 5} The $I$-pair of [1] is just the $I$-pair constructed from
the $\xvr$-structure on $\so(3)$ from remark 4.
\endremark

\subhead 3. Lie $\xvr$-algebras, Lie $R\varrho$-algebras and Lie
$\Theta\varrho$-algebras\endsubhead
First, let us describe a relation between Lie $\xvr$-algebras and Lie
$R\varrho$-algebras of the article [7].

\definition{Definition 5A {\rm [7]}} The {\it Lie $R\varrho$-algebra\/} is a
triple $(\frak g,R,\varrho)$, where $\frak g$ is the Lie algebra with the
bracket $[\cdot,\cdot]$ and $R$, $\varrho$ are two operators in it such that
the following two identities
$$
\aligned
\varrho[x,y]_{\varrho}&=[\varrho x,\varrho y],\\
R[x,y]_{\varrho}+\varrho[x,y]_R&=[Rx,\varrho y]+[\varrho x,Ry].
\endaligned
$$
holds for all $x$ and $y$ from $\frak g$. Here
$$\aligned
[x,y]_R&=[Rx,y]+[x,Ry]-R[x,y],\\
[x,y]_{\varrho}&=[\varrho x,y]+[x,\varrho y]-\varrho[x,y]+[Rx,Ry]-R[x,y]_R.
\endaligned
$$
A Lie $R\varrho$-algebra $(\frak g,R,\varrho)$ is called {\it regular\/}
if the identity
$$R[x,y]_R=2([\varrho x,y]+[x,\varrho y])$$
holds.
\enddefinition

Lie $R\varrho$-algebras are deeply related to certain quadratic bunches of
Lie algebras [7]. In particular, the bracket $[\cdot,\cdot]_\varrho$ is a Lie
one (i.e. obeys the Jacobi identity).

\proclaim{Proposition 2} Let $(\frak g,R,\varrho)$ be a regular Lie
$R\varrho$-algebra with commuting $R$ and $\varrho$. Let us define the
sequence of operators $R_n$ as
$$R_{n+1}=RR_n-\varrho R_{n-1},\quad R_0=1,\quad R_1=R.$$
The triples $(\frak g,R_n,\varrho^n)$ are regular Lie $R\varrho$-algebras.
\endproclaim

\proclaim{Proposition 3} Let $(\frak g,\xi,\varrho)$ be a Lie $\xvr$-algebra.
Let us define the sequences of operators $\xi_n$ and $R_n$ as
$$\aligned
R_{n+1}&=\Theta R_n-\varrho^2R_{n-1},\quad R_0=1,\quad R_1=\Theta;\\
\xi_{n+1}&=\Theta\xi_n+\varrho^2\xi_{n-1},\quad \xi_0=\xi,\quad
\xi_1=(\Theta+\varrho)\xi,\endaligned$$
where $\Theta=\varrho+\xi^2$. The triples $(\frak g,R_n,\varrho^{2n})$ are
regular Lie $R\varrho$-algebras, whereas the triples $(\frak g,\xi_n,
\varrho^{2n+1})$ are Lie $\xvr$-algebras.
\endproclaim

In particular, the triple $(\frak g,\Theta,\varrho^2)$ is a Lie
$R\varrho$-algebra for any Lie $\xvr$-algebra $(\frak g,\xi,\varrho)$.
So algebras from the least class may be characterized as ``square roots''
of algebras from the first one. Note once more that Lie $R\varrho$-algebras
admit an interpretation in terms of quadratic bunches of Lie algebras [7],
however, an ana\-lo\-gous interpretation for the Lie $\xvr$-algebras is not
known.

\definition{Definition 5B} The {\it Lie $\Theta\varrho$-algebra\/} is a
triple $(\frak g,\Theta,\varrho)$, where $\frak g$ is the Lie algebra with
the bracket $[\cdot,\cdot]$ and $\Theta$, $\varrho$ is two commuting
operators in it such that the following identites
$$\aligned
\varrho[x,y]_\Theta=[\varrho x,\varrho y],\quad
&\Theta[x,y]_\Theta=[\varrho^2 x,y]+[x,\varrho^2 y],\\
[\Theta x,\varrho y]+[\varrho x,\Theta y]=\varrho[\varrho x,y]
&+\varrho[x,\varrho y]+[x,\varrho y]_\Theta+[\varrho x,y]_\Theta,\\
[\varrho x,\varrho y]_\Theta=\varrho\bigl([\Theta x,\Theta y]&-
[\varrho^2 x,y]-[x,\varrho^2 y]-\varrho^2[x,y]\bigr)\\
\endaligned$$
holds for all $x$ and $y$ from $\frak g$. Here
$$[x,y]_\Theta=\tfrac12\bigl([\Theta x,y]+[x,\Theta y]-\Theta[x,y]\bigr).$$
A Lie $\Theta\varrho$-algebra is called {\it special\/} if the identity
$$\bigl([\ad x,\Theta]-\ad\Theta x\bigr)\cdot\varrho+
2\varrho\cdot\ad(\varrho x)=0$$
holds.
\enddefinition

\remark{Example 3} Let $\frak A$ be an associative algebra, $\frak A_{[\cdot,
\cdot]}$ be its commutator algebra, $\Theta x=q^2x+xq^2$, $\varrho x=qxq$
($q\in\frak A$). The triple $(\frak A_{[\cdot,\cdot]},\Theta,\varrho)$
is a special Lie $\Theta\varrho$-algebra.
\endremark

\remark{Remark 5} The bracket $[\cdot,\cdot]_\Theta$ obeys the Jacobi
identity.
\endremark

\proclaim{Proposition 4} Let $(\frak g,\Theta,\varrho)$ be a Lie
$\Theta\varrho$-algebra. Let us define the sequence of operators $\Theta_n$ as
$$\Theta_{n+1}=\Theta\Theta_n-\varrho^2\Theta_{n-1},\quad \Theta_0=1,\quad
\Theta_1=\Theta.$$
The triples $(\frak g,\Theta_n,\varrho^{2n})$ are Lie $\Theta\varrho$-algebras.
\endproclaim

Let us now relate the Lie $\Theta\varrho$-algebras to Lie $\xvr$-algebras
and Lie $R\varrho$-algebras.

\proclaim{Proposition 5}

{\bf A.} Let $(\frak g,\xi,\varrho)$ be a Lie $\xvr$-algebra. The triple
$(\frak g,\Theta,\varrho)$ ($\Theta\!=\!2\varrho\!+\!\xi^2$) is the Lie
$\Theta\varrho$-algebra.

{\bf B.} Let $(\frak g,R,\varrho)$ be a regular Lie $R\varrho$-algebra with
commuting $R$ and $\varrho$. The triple $(\frak g,\Theta,\varrho)$
($\Theta\!=\!R^2\!-\!2\varrho$) is the Lie $\Theta\varrho$-algebra.
\endproclaim

It is rather interesting to describe the classes of Lie $\xvr$-algebras and
regular $R\varrho$-algebras, to which the special $\Theta\varrho$-algebras
are related.

\Refs
\roster
\item"[1]" Juriev D., On the nonHamiltonian interaction of two rotators.
E-print: dg-ga/9409004.
\item"[2]" Juriev D., An excursus into the inverse problem of
representation theory [in Russian]. E-print: mp\_arc/96-477.
\item"[3]" Maltsev A.I., Algebraic systems. Moscow, 1970.
\item"[4]" Koecher M., Jordan algebras and their applications, 1962.
\item"[5]" Koecher M., An elementary approach to bounded symmetric domains,
Houston, 1969.
\item"[6]" Juriev D., Topics in isotopic pairs and their representations.
III. Bunches of Lie algebras and modified classical Yang-Baxter equation.
E-print: q-alg/9708027.
\item"[7]" Juriev D., Topics in hidden symmetries. VI. E-print: q-alg/9708028.
\endroster
\endRefs
\enddocument